\documentclass[1p,review,authoryear]{amsart}

\usepackage[latin1]{inputenc}
\usepackage[T1]{fontenc}
\usepackage{amsmath,amssymb,amsfonts}
\usepackage{graphicx}
\usepackage{array}
\usepackage{float}
\makeatletter
\def\ps@pprintTitle{%
  \let\@oddhead\@empty
  \let\@evenhead\@empty
  \let\@oddfoot\@empty
  \let\@evenfoot\@oddfoot
}

\begin{document}

\title{Alternative to the Romberg Method of Estimating the Definite Integral}

\begin{center}
\author{M. Brandon Youngberg}
\end{center}

\begin{abstract}
\noindent Using elementary methods, we define and derive a particular weighted average of the trapezoidal and composite trapezoidal rules and show that this approximation, as well as its composite, is straightforward in computation.  This approximation and its composite, in their general forms, are shown to have predictable error patterns; thus, an extrapolation method can be used to increase the accuracy.  We then derive the necessary weights to use an extrapolation method to reduce error and converge more quickly than Romberg integration by allowing for improved accuracy with fewer necessary subintervals.  The procedure necessary to implement this alternative method is then carefully described, followed by two examples.\\
\end{abstract}

\maketitle

\section{Introduction}

\noindent A variety of numerical integration methods exist with diverse rates of convergence, codability, and ease of computation, and with varying requirements, both of the function to be integrated as well as the information required for approximation.  In this paper, we will examine a new method of integration similar in application and scope to the Romberg method of integration, but with relatively faster rates of convergence because of fewer information requirements. \\

\noindent After providing motivation for its construction, we introduce a form of the initial approximation, denoted $A$, with emphasis on the ease of this approximation's calculation.  We derive $A$, a particular combination of the trapezoidal and composite trapezoidal rules, and then the remaining error of $A$, denoted $\alpha$.  This error is derived by weighting the errors of the trapezoidal rule and composite trapezoidal rule by the same weights used in the original approximation.\\

\noindent Having $A$, we derive its composite, $A_m$, and determine its error, after which it becomes apparent that an extrapolation technique could be applied.  We then derive the form of $\Omega$, the extrapolation factor necessary to eliminate the first error term between $A$ and $A_m$.  We then find the general form of these possible extrapolation factors that are necessary to eliminate that first error term for all feasible combinations of $A$ and $A_m$.  We perform the first extrapolation and then we find the next few extrapolation factors necessary to eliminate further error terms.\\

\noindent Afterwards, we provide in detail the procedure necessary to appropriately use the extrapolation factors and approximations to arrive at the most accurate approximation, before finally including two examples on this technique's implementation.

\section{Motivation}
\noindent Consider the following graph, Figure 1, of a monotonically increasing function.  Often such a graph is used to introduce numerical integration concepts, particularly when the underestimate and overestimate rectangles are included, and the stepsizes, denoted $h$, are consistent as they are here.  As is often the case, we define the stepsize as $h=(b-a)/n$ where $b-a$ is the entire width of the integral, and $n$ is the number of subintervals being used.  Also, we define the sum of these underestimate rectangles as $U_1$, and the sum of these overestimate rectangles as $O_1$.\\

\begin{figure}[H]
\centering
\includegraphics[scale=0.65]{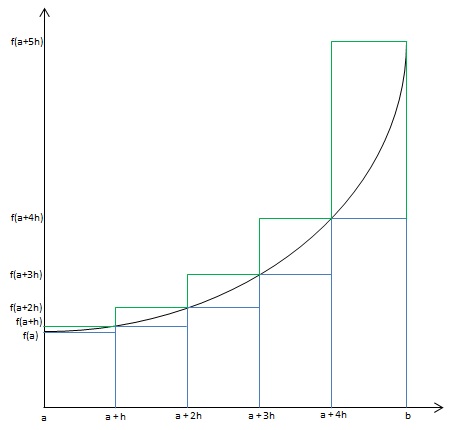}
\caption{}
\end{figure}

\noindent Now consider if we were to segregate the differences between the underestimate and overestimate rectangles and to stack these differences.  Because the stepsizes are equal, we know that the width of this stack will have the value $h$.  Additionally, if we include the height of the lowest value for $f(x)$, which here is at $f(a)$, we have a column that is exactly the same height and width as the largest overestimate rectangle.  This concept is depicted in Figure 2 below.\\

\begin{figure}[H]
\centering
\includegraphics[scale=0.65]{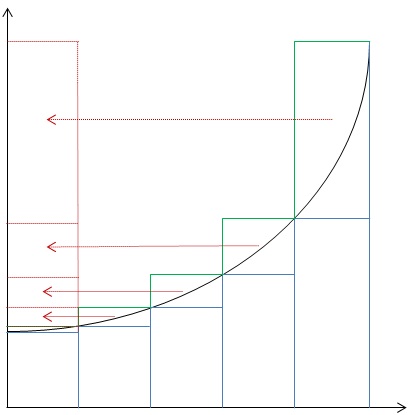}
\caption{}
\end{figure}

\noindent Now consider if we were to divide the width of this new rectangle by the same number of underestimate or overestimate rectangles that were in the original graph.  We could then adjust the original function appropriately to fit in this new rectangle.  The original values for $f(x)$ would be in their original position, but now $h = (b-a)/n^2$. This concept is depicted in the following graph, Figure 3.  \\

\begin{figure}[H]
\centering
\includegraphics[scale=0.65]{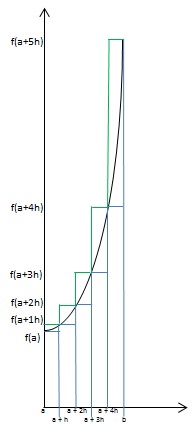}
\caption{}
\end{figure}

\noindent Intuitively, here the underestimate rectangles could be regarded as underestimates of the error of the original underestimates.  Adding these underestimates to the original underestimates provides a more accurate estimate.  Similarly, these new overestimates are overestimates of the error of the original overestimates.  We must subtract these from the original overestimates to provide a more accurate estimate.  Now consider if we were to repeat this process an infinite number of times, each iteration more narrow than the previous.  We could then sum these iterations for a more accurate underestimate and overestimate.  This sum for the underestimate will be referred to as $U_{\infty}$, and for the overestimate, $O_{\infty}$.  Our approximation will be the average of $U_{\infty}$ and $O_{\infty}$.  Additionally, it will be shown in the next section that, strangely, for this approximation to be valid, it is not necessary for the function to be monotonic at all!

\section{Derivation of the Initial Approximation and Initial Error}
\noindent For $i \in \{0,1,\ldots, n\}$, set $x_i = a + ih$. Define $U_1$ to be the first underestimate, so that $U_1 = h\sum\limits_{i = 0}^{n-1} f(x_i)$, and $O_1$ to be the first overestimate, so that $O_1 = h \sum\limits_{i = 1}^{n} f(x_i)$. For future reference note that if $\Sigma = \sum\limits_{i = 1}^{n-1} f(x_i)$, then

\begin{align*}
U_1 + O_1 & = h(f(a) + 2\Sigma + f(b)), \\
U_1 - O_1 & = h(f(a) - f(b)).
\end{align*}

\noindent Let $U_{\infty}$ and $O_{\infty}$ denote the limit of the underestimates and overestimates, respectively. Finally, set $H = hf(a)$. We then have

\begin{align}
U_{\infty} & = U_1 + \dfrac{U_1}{n-1} - \dfrac{n}{n-1}H, \\
O_{\infty} & = O_1 - \dfrac{O_1}{n+1} + \dfrac{n}{n+1}H.
\end{align}

\noindent The average of $(1)$ and $(2)$ is then

\begin{align}
& A = \dfrac{U_{\infty} + O_{\infty}}{2}, \nonumber \\
& = \dfrac{1}{2}\left[U_1 + O_1 + \dfrac{U_1}{n-1} - \dfrac{O_1}{n+1} + \left(\dfrac{n}{n+1} - \dfrac{n}{n-1}\right)H\right], \nonumber \\
& = \dfrac{1}{2}\left[\dfrac{(n^2-1)(U_1 + O_1) + (n+1)U_1 - (n-1)O_1}{n^2 - 1} - \dfrac{2n}{n^2 - 1}H\right], \nonumber \\
& = \dfrac{1}{2}\left[\dfrac{n^2(U_1 + O_1) + n(U_1 + O_1) - 2nH}{n^2-1}\right], \nonumber \\
& = \dfrac{1}{2}\left[\dfrac{hn^2(f(a) + 2\Sigma + f(b)) + hn(f(a) - f(b)) - 2hnf(a)}{n^2 - 1}\right], \nonumber \\
& = \dfrac{n^2\left[\dfrac{h}{2}(f(a) +2\Sigma + f(b)\right] - \dfrac{hn}{2}(f(a) + f(b))}{n^2-1}.
\end{align}

\noindent Now note that $\frac{h}{2}(f(a) +2\Sigma + f(b)) = \frac{b-a}{2n}(f(a) +2\Sigma + f(b))$ is the estimate of $\int_a^b f(x)\ dx$ coming from the composite trapezoidal rule and that $\frac{hn}{2}(f(a) + f(b)) = \frac{b-a}{2}(f(a) + f(b))$ is the estimate of $\int_a^b f(x)\ dx$ coming from the ordinary trapezoidal rule. \\

\noindent To find our error, we can use Taylor expansion to derive the trapezoidal rule.  We expand around $a$ and $b$, and then we combine the expressions to obtain

\begin{align*}
\displaystyle\int_a^b f(x)\ dx = \dfrac{h}{2}[f(a)+f(b)]+\dfrac{h^2}{2(2!)}[f^{'}(a)-f^{'}(b)]+\dfrac{h^3}{2(3!)}[f^{''}(a)+f^{''}(b)]+O(h^4),
\end{align*}

\noindent where $h=b-a$.  We can then follow a well-known procedure to further use Taylor expansion and remove odd powers of $h$, which gives us the Euler-Maclaurin summation formula shown here:

\begin{align}
\displaystyle\int_a^b f(x)\ dx = \dfrac{h}{2}[f(a)+f(b)]+\dfrac{h^2}{6(2!)}[f^{'}(a)-f^{'}(b)]-\dfrac{h^4}{(6!)}[f^{'''}(a)-f^{'''}(b)]+O(h^6).
\end{align}

\noindent Subtracting the trapezoidal rule from $(4)$ results in the trapezoidal rule error.  Denote the first term of this error as $E_{T1}$, such that $E_{T1} = \frac{h^2}{6(2!)}[f^{'}(a)-f^{'}(b)]$.  Further discussion of the derivation the trapezoidal rule, its error, and the Euler-Maclaurin summation formula can be found in [2, 4, 5]. \\

\noindent Now when we construct a composite, intermediate points cancel.  To demonstrate this, consider the construction of a composite trapezoidal rule with two parts:

\begin{align*}
\displaystyle\int_{x_1}^{x_3} f(x)\ dx = \dfrac{h}{2}[f(x_1)+f(x_2)]+\dfrac{h}{2}[f(x_2)+f(x_3)] \nonumber \\ +\dfrac{h^2}{6(2!)}[f^{'}(x_1)-f^{'}(x_2)]+\dfrac{h^2}{6(2!)}[f^{'}(x_2)-f^{'}(x_3)]+O(h^4).
\end{align*}

\noindent Notice that in the last two terms, the $f({x_2})$ terms cancel.  This occurs for each set of derivatives as they are always expressed as differences using this derivation method.  Additionally, similar to as how the trapezoidal rule uses the entire width of the integral as its stepsize, the composite trapezoidal rule uses the entire width of each subinterval as its stepsize.  If $w$ is the width of each of these subintervals, then in the composite trapezoidal rule, $h=nw/n$, but $nw = b-a$.  So if in the trapezoidal rule error we replace $h=(b-a)$ with $h=(b-a)/n$, we are left with the composite trapezoidal rule error.  Then, this composite trapezoidal rule error is written as

\begin{align*}
\dfrac{(b-a)^2}{6(2!)n^2}[f^{'}(a)-f^{'}(b)]-\dfrac{(b-a)^4}{(6!)n^4}[f^{'''}(a)-f^{'''}(b)]+O(h^6).
\end{align*}

\noindent Denote the first term of this error as $E_{C1}$.  Then, notice that $E_{C1}$, when multiplied by $n^2$, is exactly the same as $E_{T1}$.  Then we can rewrite expression $(3)$ to eliminate these first error terms and to derive the remaining error, denoted as $\alpha$, as follows:

\begin{align*}
A = \dfrac{n^2\left[\dfrac{h}{2}(f(a) +2\Sigma + f(b)\right] - \dfrac{hn}{2}(f(a) + f(b))}{n^2-1}, \\
= \dfrac{n^2\left[\displaystyle\int_a^b f(x)\ dx + E_C\right] - \left[\displaystyle\int_a^b f(x)\ dx + E_T\right]}{n^2-1},
\end{align*}

\begin{align*}
= \dfrac{n^2\left[\displaystyle\int_a^b f(x)\ dx + \dfrac{(b-a)^4}{6!n^4}[f^{'''}(a) - f^{'''}(b)]+ O(h^6)\right] - \dfrac{n^2(b-a)^2}{2(3!)n^2}[f^{'}(a) - f^{'}(b)]}{n^2-1} \\
+ \dfrac{\dfrac{n^2(b-a)^2}{2(3!)n^2}[f^{'}(a) - f^{'}(b)] - \left[\displaystyle\int_a^b f(x)\ dx + \dfrac{n^4(b-a)^4}{6!n^4}[f^{'''}(a) - f^{'''}(b)]+ O((b-a)^6)\right]}{n^2-1}, \\
= \dfrac{(n^2-1)\left[\displaystyle\int_a^b f(x)\ dx\right] - n^2(n^2-1)\left[\dfrac{(b - a)^4}{6!n^4}[f^{'''}(a) - f{'''}(b)]\right]}{n^2-1}+O\left(g(h)\right).
\end{align*}

\noindent Then, $A = \displaystyle\int_a^b f(x)\ dx - \dfrac{(b - a)^4}{6!n^2}[f^{'''}(a)-f^{'''}(b)] +O\left(g(h)\right)$, and 

\begin{align}
& \alpha = - \dfrac{(b - a)^4}{6!n^2}[f^{'''}(a)-f^{'''}(b)] +O\left(g(h)\right).
\end{align}

\noindent where 

\begin{align*}
g(h) = c\left(\dfrac{n^2h^6-(b-a)^6}{n^2-1}\right) = c\left(\dfrac{-n^2(b-a)^6-(b-a)^6}{n^4}\right),
\end{align*}

\noindent and where $c$ is the arbitrary constant.  It will be shown in the following section that $g(h)$ will be able to be further simplified quite nicely, and it has been purposefully left in this form to more easily facilitate this simplification.

\section{Derivation of the Composite and the Extrapolation Factors}

\noindent Denote the composite approximation as $A_m$, where $m$ is the number of subintervals used in each part of the composite.  We can clarify the relationship between the original approximation's error, $A$, and that of the composite approximation, $A_m$, by first observing that the original approximation's error is such that

\begin{align*}
\alpha = - \dfrac{(b - a)^4}{6!n^2}[f^{'''}(a)-f^{'''}(b)]+O\left(g(h)\right),\\
= - \dfrac{n^2(b - a)^4}{6!n^4}[f^{'''}(a)-f^{'''}(b)]+O\left(g(h)\right),\\
= - \dfrac{n^2h^4}{6!}[f^{'''}(a)-f^{'''}(b)]+O\left(g(h)\right).
\end{align*}

\noindent Consider, again, if we were to construct a composite with $k$ parts of equal width $w$.  Then, since each $w$ is equivalent, $kw = b - a$ and $km = n$.  This allows us to now define our stepsize for the entire integral as $h = w/m = kw/km = (b-a)/n$.  Then, each part of our composite would have the following error, denoted $\alpha_m$, where the subscript $m$ corresponds to the number of subintervals used in each part of the composite:

\begin{align*}
\alpha_m = - \dfrac{m^2w^4}{6!m^4}[f^{'''}(a)-f^{'''}(b)]+O(j(h)), 
\end{align*}

\noindent where

\begin{align*}
j(h)= c\left(\dfrac{m^2h^6-(b-a)^6}{n^2-1}\right) = c\left(\dfrac{-m^2(b-a)^6-(b-a)^6}{n^4}\right).
\end{align*}

\noindent Now the composite error for the entire integral becomes

\begin{align}
\alpha_m = - \dfrac{m^2(b - a)^4}{6!n^4}[f^{'''}(a)-f^{'''}(b)]+O\left(j(h)\right).
\end{align}

\noindent It is now apparent that an extrapolation is possible since we will have two approximations with comparable errors, the original and a composite.  Then we must first solve for the first weight, or extrapolation factor, used in our extrapolation process.  We will use $\Omega$ to denote these extrapolation factors, with subscripts to carefully track for which extrapolation each $\Omega$ is used.  Calculation of $\Omega_{n,m}$, where the subscript $n$ denotes the total number of subintervals in the entire integral and where $m$ denotes the number of subintervals used in the construction of the composite, is as follows:

\begin{align}
& \alpha = \Omega_{n,m} \alpha_m, \nonumber \\
& - \dfrac{(b - a)^4}{6!n^2}[f^{'''}(a)-f^{'''}(b)]= \Omega_{n,m} \left[-\dfrac{m^2(b - a)^4}{6!n^4}[f^{'''}(a)-f^{'''}(b)]\right], \nonumber \\
& \Omega_{n,m} = \left(\dfrac{n}{m}\right)^2.
\end{align}

\noindent  This holds providing $m$ divides $n$ evenly and $2 \le m \le n/2$, as a larger $m$ disallows construction of composites with parts of equal numbers of subintervals.  From this point forward, we shall denote $A$ as $A_n$ and its error as $\alpha$ as $\alpha_n$, because we can regard $A_n$ as a composite constructed of one part with $n$ subintervals.  \\

\noindent We can now reduce the error from $O(n^2h^4)$ to $O(n^2h^6)$ provided we use $A_n$ and any feasible $A_m$, providing the composite is constructed with parts of equal numbers of subintervals.  This straightforward procedure is conducted as follows:

\begin{align*}
& \displaystyle\int_a^b f(x)\ dx = \dfrac{\Omega_{n,m} A_m - A_n}{\Omega_{n,m} - 1},\\
& = \dfrac{\left(\dfrac{n}{m}\right)^2\left[\displaystyle\int_a^b f(x)\ dx - \alpha_m\right] - \left[\displaystyle\int_a^b f(x)\ dx - \alpha_n  \right]}{\left(\dfrac{n}{m}\right)^2-1},\\
& = \dfrac{\left(\dfrac{n}{m}\right)^2\left[\displaystyle\int_a^b f(x)\ dx + \dfrac{m^2(b - a)^4}{6!n^4}[f^{'''}(a)-f^{'''}(b)]+O\left(j(h)\right)\right]}{\left(\dfrac{n}{m}\right)^2-1}\\
& - \dfrac{\left[\displaystyle\int_a^b f(x)\ dx + \dfrac{(b - a)^4}{6!n^2}[f^{'''}(a)-f^{'''}(b)]+O\left(g(h)\right) \right]}{\left(\dfrac{n}{m}\right)^2-1},\\
& = \displaystyle\int_a^b f(x)\ dx + O\left(\dfrac{\left(\dfrac{n}{m}\right)^2j(h)-g(h)}{\left(\dfrac{n}{m}\right)^2-1}\right),\\
& = \displaystyle\int_a^b f(x)\ dx + O\left(c\left(\dfrac{(-m^2n^2-n^2)(b-a)^6+(m^2n^2+m^2)(b-a)^6}{n^4(n^2-m2)}\right)\right),\\
& = \displaystyle\int_a^b f(x)\ dx + O(n^2h^6).
\end{align*}

\noindent Additionally, note that further error terms are described generally by $O(n^2h^\rho)$ where $\rho$ is even.  This is important to recognize, as an increase in the number of points available will not necessarily increase the accuracy by an expected amount.  For instance, should the first error term be simply $O(h^6)$, then doubling the number of points will increase the accuracy by a factor of 64.  However, doubling the number of available points when the first error term is described as $O(n^2h^6)$ will only initially increase the accuracy by a factor of 16.  Incidentally, doubling the number of available points would also allow for the use of an additional composite, which would ultimately increase the accuracy, as will be demonstrated in the following sections.\\
\\
\noindent At this point, it is interesting to note that, after some algebra, we can rewrite our approximation $A_n$ in the following form:

\begin{align}
A_n = \dfrac{h}{2}\left[\dfrac{nf(a)}{n+1}+\dfrac{2n^2\Sigma}{n^2-1}+\dfrac{nf(b)}{n+1}\right].
\end{align}

\noindent Then, when $n=2$, $(8)$ is exactly Simpson's rule with $(5)$ and $(6)$ exactly Simpson's Rule error and composite Simpson's Rule error, respectively.  Additionally, when $n=3$, the same is true concerning Simpson's 3/8 Rule and errors as well.  However, beyond the case where $n=3$, $(8)$ ceases to follow further Newton-Cotes formulas exactly.  Further discussion of the Newton-Cotes formulas and errors, including Simpson's Rule and Simpson's 3/8 Rule can be found in [1, 6].\\

\noindent The extrapolation factors necessary to remove the remaining error terms are very simple to compute.  Consider two separate composites, one comprised of parts of $x$ subintervals, denoted $A_x$, and one comprised of parts of $y$ subintervals, denoted $A_y$.  Then let the remaining error terms of these composites to be $\alpha_x$ and $\alpha_y$, respectively.  As shown in (7), we will have two separate extrapolation factors, $\Omega_{n,x}$ and $\Omega_{n,y}$.  Then we know that

\begin{align}
\alpha_n = \Omega_{n,x} \alpha_x, \nonumber
\end{align}

\noindent and also that
\begin{align}
\alpha_n = \Omega_{n,y} \alpha_y. \nonumber
\end{align}

\noindent This allows us to then compute

\begin{align}
\alpha_y = \dfrac{\Omega_{n,x}}{\Omega_{n,y}} \alpha_x = \left(\dfrac{y}{x}\right)^2\alpha_x.
\end{align}

\noindent We can now extend the use of the extrapolation factors' subscripts to denote combinations of previously-used extrapolation factors from which they result, such as:

\begin{align}
\Omega_{y,x} = \dfrac{\Omega_{n,x}}{\Omega_{n,y}}=\left(\dfrac{n}{x}\right)^2\left(\dfrac{y}{n}\right)^2=\left(\dfrac{y}{x}\right)^2. \nonumber
\end{align}

\noindent Now consider if we had a third composite composed of parts of $z$ subintervals.  The associated error would then be denoted $\alpha_z$.  We could then calculate $\alpha_z$, similarly as to $\alpha_x$ in (9), giving us

\begin{align}
\alpha_z = \Omega_{z,y}\alpha_y = \dfrac{\Omega_{n,y}}{\Omega_{n,z}} \alpha_y = \left(\dfrac{z}{y}\right)^2\alpha_y. \nonumber
\end{align}

\noindent Furthermore,

\begin{align}
\alpha_z = \left(\dfrac{z}{y}\right)^2\alpha_y = \left(\dfrac{z}{y}\right)^2\left(\dfrac{y}{x}\right)^2\alpha_x = \Omega_{z,y}\Omega_{y,x} \alpha_y = \left(\dfrac{z}{x}\right)^2\alpha_x = \Omega_{z,x}\alpha_x. \nonumber
\end{align}

\noindent Then generally,

\begin{align*}
\alpha_z = \left(\dfrac{z}{y}\right)^2\alpha_y = \left(\dfrac{z}{y}\right)^2\left(\dfrac{y}{x}\right)^2... \left(\dfrac{c}{b}\right)^2\left(\dfrac{b}{a}\right)^2\alpha_a = \Omega_{z,y}\Omega_{y,x} ... \Omega_{c,b}\Omega_{b,a} = \left(\dfrac{z}{a}\right)^2\alpha_a = \Omega_{z,a}\alpha_a. \nonumber
\end{align*}

\noindent Therefore, any $\alpha$ can be combined with an appropriate extrapolation factor, $\Omega$, where the $\Omega$ used is simply a combination of previously-used $\Omega$, in such a way as to equate this $\alpha$ with $\alpha_z$.  If we allow $z = n$, where $n$ is the total number of subintervals available in the entire integral, it becomes clear that we have, therefore, found the general form of all of the extrapolation factors necessary to reduce the error to the smallest value that the number of available composites will allow.  Additionally, the general extrapolation factor $\Omega_{z,a}$ is sufficient to convert all remaining error terms of our new approximation from $O(a^2h^{\rho})$ to  $O(z^2h^{\rho})$.

\section{Procedure}
\noindent Now that we have specified how to compute $A_n$, $A_m$, and the extrapolation factor $\Omega$, we use the following procedure to produce the most accurate estimate available.  Throughout this procedure, we must keep careful track of the subscripts of our extrapolation factors, as well as how we apply these factors in the extrapolation.  Failure to carefully note the number of subintervals used in the construction of a composite, and how its associated extrapolation factor is used, may cause error in calculation of further extrapolations.  Table 1 at the end of this section is provided as an aid for following the procedure.\\

\noindent The first step is to compute $A_n$, using the average of $(1)$ and $(2)$ as described in section 3, as well as $A_m$ for all available composites.  The composites must be constructed with parts comprised of an equal number of subintervals.  Though these composites, as well as the initial approximation, will likely be different, they all have their first error term of $O(m^2h^4)$, where $m$ is the number of subintervals used in each composite.\\

\noindent Suppose there are four separate available composites, which we will denote as $A_a$, $A_b$, $A_c$, and $A_d$,  where the number of subintervals used in each part of each composite is $a$, $b$, $c$, and $d$ respectively. As previously shown in $(7)$, the $\Omega$ for the extrapolation between $A_n$ and $A_a$ will be $(n/a)^2$, between $A_a$ and $A_b$ will be $(a/b)^2$, between $A_b$ and $A_c$ will be $(b/c)^2$, and between $A_c$ and $A_d$ will be $(c/d)^2$.\\

\noindent Perform the first extrapolation for each composite.  Remember that care must be used to track the extrapolation factors for use in further extrapolations and computations of $\Omega$.  The numerator of the extrapolation factor corresponds to the approximation being subtracted in the extrapolation process.  For instance, the first extrapolation will be conducted using the original approximation, $A_n$ and a composite as follows:

\begin{align*}
A_{n,a} = \dfrac{\Omega_{n,a} A_a - A_n}{\Omega_{n,a}- 1}= \dfrac{\left(\dfrac{n}{a}\right)^2 A_a - A_n}{\left(\dfrac{n}{a}\right)^2-1}.
\end{align*}

\noindent Notice again that the numerator of the extrapolation factor, in this case $n$, corresponds with the approximation being subtracted.  The next approximation would be

\begin{align*}
A_{a,b} = \dfrac{\Omega_{a,b} A_b - A_a}{\Omega_{a,b} - 1} = \dfrac{\left(\dfrac{a}{b}\right)^2 A_b - A_a}{\left(\dfrac{a}{b}\right)^2-1}.
\end{align*}

\noindent Continuing in this manner, we should have several new approximations, one for each extrapolation we have performed so far.  Denote these new approximations $A_{n,a}$, $A_{a,b}$, $A_{b,c}$, and $A_{c,d}$.  Now we compute new extrapolation factors.  As was previously shown in $(9)$, the extrapolation factor to reduce the error using approximations $A_{n,a}$ and $A_{a,b}$ is $\Omega_{n,b}=(n/a)^2(a/b)^2=(n/b)^2$.  Now perform the next set of extrapolations as shown here with $A_{n,a}$ and $A_{a,b}$:

\begin{align*}
A_{n,b} = \dfrac{\Omega_{n,b} A_{a,b} - A_{n,a}}{\Omega_{n,b} - 1}
=\dfrac{\left(\dfrac{n}{b}\right)^2 A_{a,b} - A_{n,a}}{\left(\dfrac{n}{b}\right)^2-1}.
\end{align*}

\noindent Again, we have new approximations, $A_{n,b}$, $A_{a,c}$, and $A_{b,d}$.  Compute the next extrapolation factors, $\Omega_{n,c}$ and $\Omega_{a,d}$.  Then compute the next set of extrapolations and obtain the next two approximations.  This procedure conducted with approximations $A_{n,b}$ and $A_{a,c}$ using $\Omega_{n,c}=(n/a)^2(a/b)^2(b/c)^2=(n/c)^2$ is done so in the following manner:

\begin{align*}
A_{n,c} = \dfrac{\Omega_{n,c} A_{a,c} - A_{n,b}}{\Omega_{n,c} - 1} 
=\dfrac{\left(\dfrac{n}{c}\right)^2 A_{a,c} - A_{n,b}}{\left(\dfrac{n}{c}\right)^2-1}.
\end{align*}

\noindent Two new approximations result, $A_{n,c}$ and $A_{a,d}$.  Calculate the next extrapolation factor, $\Omega_{n,d} =(n/a)^2(a/b)^2(b/c)^2(c/d)^2=(n/d)^2$.  Then,

\begin{align*}
A_{n,d} = \dfrac{\Omega_{n,d} A_{a,d} - A_{n,c}}{\Omega_{n,d} - 1}
=\dfrac{\left(\dfrac{n}{d}\right)^2 A_{a,d} - A_{n,c}}{\left(\dfrac{n}{d}\right)^2-1}.
\end{align*}

\noindent This provides the highest accuracy estimation we are able to achieve given the original approximation, $A_n$ and four composites.  Note that this final approximation will have its first error term as $O(n^2h^{12})$, because all error terms of the form $O(m^2h^{12})$ have been replaced with the error corresponding to the composite composed of parts of $n$ subintervals, the original approximation.\\

\noindent The following table helps to keep track of the extrapolation procedure and to clarify the manner in which this procedure causes convergence.  Additionally, when implementing this procedure it is tempting to further reduce the values for $\Omega$ or to express them in decimal form, but we have left them reduced only to the level used in this procedure, as doing so may help to reduce confusion when constructing further $\Omega$.

\begin{table}[H]
\centering
\caption{}
\scalebox{0.65}{
\begin{tabular}{r c c c c c c c c c c}

Error & $O(n^2h^4)$ & &  $O(n^2h^6)$ & & $O(n^2h^8)$ & & $O(n^2h^{10})$ & & $O(n^2h^{12})$ \\ [0.5ex]
\hline
& $A_n$ & $\longrightarrow$ & $A_{n,a}= \dfrac{\left(\dfrac{n}{a}\right)^2 A_a - A_n}{\left(\dfrac{n}{a}\right)^2-1}$ & $\longrightarrow$ & $A_{n,b} = \dfrac{\left(\dfrac{n}{b}\right)^2 A_{a,b} - A_{n,a}}{\left(\dfrac{n}{b}\right)^2-1}$ & $\longrightarrow$ & $A_{n,c} = \dfrac{\left(\dfrac{n}{c}\right)^2 A_{a,c} - A_{n,b}}{\left(\dfrac{n}{c}\right)^2-1}$ & $\longrightarrow$ & $A_{n,d} = \dfrac{\left(\dfrac{n}{d}\right)^2 A_{a,d} - A_{n,c}}{\left(\dfrac{n}{d}\right)^2 - 1}$ & \\

& $A_a$ & $\longrightarrow \nearrow$ & $A_{a,b}= \dfrac{\left(\dfrac{a}{b}\right)^2 A_b - A_a}{\left(\dfrac{a}{b}\right)^2-1}$ & $\longrightarrow \nearrow$ & $A_{a,c} = \dfrac{\left(\dfrac{a}{c}\right)^2 A_{b,c} - A_{a,b}}{\left(\dfrac{a}{c}\right)^2-1}$ & $\longrightarrow \nearrow$ & $A_{a,d} = \dfrac{\left(\dfrac{a}{d}\right)^2 A_{b,d} - A_{a,c}}{\left(\dfrac{d}{b}\right)^2-1}$ & $\nearrow$ & \\

& $A_b$ & $\longrightarrow \nearrow$ & $A_{b,c}= \dfrac{\left(\dfrac{b}{c}\right)^2 A_c - A_b}{\left(\dfrac{b}{c}\right)^2-1}$ & $\longrightarrow \nearrow$ & $A_{b,d} = \dfrac{\left(\dfrac{b}{d}\right)^2 A_{c,d} - A_{b,c}}{\left(\dfrac{b}{d}\right)^2-1}$ & $ \nearrow$ & \\

& $A_c$ & $\longrightarrow \nearrow$ & $A_{c,d}= \dfrac{\left(\dfrac{c}{d}\right)^2 A_d - A_c}{\left(\dfrac{c}{d}\right)^2-1}$ & $\nearrow$ & \vspace{0.25in} \\
& $A_d$ & $\nearrow$ \\ [1ex]
\hline
\end{tabular}}
\label{table:nonlin}\\
\end{table}

\vspace{0.25in}
\section{Examples}

\noindent Example 1:\\

\noindent This example was chosen because we know with the use of extrapolation we can arrive at an error term that is a multiple of derivatives that have the same constant value.  This will ensure that our solution is exact.  Additionally, $h = 1$ which simplifies some of our calculations.\\

\noindent Consider $f(x) = x^7 -2x + 10$.  Then $\int_0^{10} f(x)\ dx = 12500000$.  Assume, however, that we do not know the function, but are instead given eleven equally spaced points from x = 0 to 10 and their corresponding values for $f(x)$.

\begin{table}[H]
\centering
\begin{tabular}{l l}
$x$ & $f(x)$ \\ [0.5ex]
\hline
0 & 10\\
1 & 9 \\
2 & 134 \\
3 & 2191 \\
4 & 16386 \\
5 & 78125 \\
6 & 279934 \\
7 & 823539 \\
8 & 2097146 \\
9 & 4782961 \\
10 & 9999990 \\ [1ex]
\hline
\end{tabular}
\label{table:1}
\end{table}

\noindent Our first step is to calculate $A_{10}$ by first calculating $U_\infty$ and $O_\infty$.  However, note that $b - a = 10 = n$, and therefore, $h = 1$.  Also, note that $H = hf(a) = 10$.  Then, as stated previously, we calculate $U_1$ as the underestimate, which is $h$ multiplied by the sum of the first 10 values of $f(x)$, which gives us $U_1 = 8080435$.  Likewise, the overestimate becomes $O_1 = 18080415$.  Then we have

\begin{align}
U_{\infty} & = U_1 + \dfrac{U_1}{n-1} - \dfrac{n}{n-1}H = 8080435 + \dfrac{8080435}{9} - \dfrac{10}{9}\cdot 10 = 8978250, \\
O_{\infty} & = O_1 - \dfrac{O_1}{n+1} + \dfrac{n}{n+1}H = 18080415 - \dfrac{18080415}{11} + \dfrac{10}{11}\cdot 10 = 16436750.
\end{align}

\noindent Averaging $(10)$ and $(11)$ then yields the initial approximation, $A_{10} = 12707500$.\\

\noindent Now we can construct $A_m$.  We have two options when $n = 10$, a composite constructed using two parts of five subintervals and a composite constructed using five parts of two subintervals.\\

\noindent Though the composite using five parts would provide a more accurate estimate, for the purposes of this example we will start with the composite using two parts.  In this case note that, for each of these parts, $m = 5$.  This does not change the value of $h$, however, as the width of each part of the composite is also $5$.  It does, however, change the value of $H$.  Because this composite is found by simply adding the results of the two separate approximations using five subintervals each, we will have a different $f(a)$ for each of these parts, and thus a separate $H$.  In the first part, $H_1 = hf(0) = 10$.  However, in the second part, $H_2 = hf(5) = 78125$.  We will also have two separate sets of overestimates and underestimates, one set for each part of which the composite is comprised.\\

\noindent The first underestimate is found using values $f(0)$ through $f(4)$, and our first overestimate is found using values $f(1)$ through $f(5)$.  Then our second underestimate is found using values $f(5)$ through $f(9)$, and our second overestimate is found by using values $f(6)$ through $f(10)$.  Then we have

\begin{align}
U_{\infty} & = U_1 + \dfrac{U_1}{n-1} - \dfrac{n}{n-1}H_1 = 18730 + \dfrac{18730}{4} - \dfrac{5}{4}\cdot 10 = 23400,\\
O_{\infty} & = O_1 - \dfrac{O_1}{n+1} + \dfrac{n}{n+1}H_1 = 96845 - \dfrac{96845}{6} + \dfrac{5}{6}\cdot 10 = 80712.5.
\end{align}

\noindent Averaging $(12)$ and $(13)$ gives us $52056.25$, which is our approximation of the area of the first half of the entire integral.

\begin{align}
U_{\infty} & = U_2 + \dfrac{U_2}{n-1} - \dfrac{n}{n-1}H_2 = 8061705 + \dfrac{8061705}{4} - \dfrac{5}{4}\cdot 78125 = 9979475,\\
O_{\infty} & = O_2 - \dfrac{O_2}{n+1} + \dfrac{n}{n+1}H_2 = 17983570 - \dfrac{17983570}{6} + \dfrac{5}{6}\cdot 78125 = 15051412.5.
\end{align}

\noindent Averaging $(14)$ and $(15)$ gives us $12515443.75$, which is our approximation of the area of the second half of the entire integral.  Adding this to the approximation of the first half gives us $A_5 = 12567500$.\\

\noindent Now that we have our initial approximation and a composite approximation, we can find the appropriate extrapolation factor $\Omega_{10,5}$.  As this is our first extrapolation, we know that the numerator of $\Omega_{10,5}$ is the total number of subintervals in the original approximation and that the denominator is the number of subintervals used in each part of the composite.  Therefore, in this case, $\Omega_{10,5} = (n/m)^2=(10/5)^2$.\\

\noindent We calculate our new approximation, $A_{10,5}$, as follows:

\begin{align*}
A_{10,5} = \dfrac{\left(\dfrac{10}{5}\right)^2A_5 - A_{10}}{\left(\dfrac{10}{5}\right)^2-1} =\dfrac{\left(\dfrac{10}{5}\right)^212567500 - 12707500}{\left(\dfrac{10}{5}\right)^2-1} = 12520833.\overline{3}.
\end{align*}

\noindent Calculating a composite using parts containing two subintervals each gives us $A_2 = 12511500$.  We then find $A_{5,2} = 12500833.\overline{3}$ using $\Omega_{5,2} = \Omega_{10,2}/\Omega_{10,5}=(10/2)^2/(10/5)^2=(5/2)^2$.  We are now able to further combine these two approximations, $A_{10,5}$ and $A_{5,2}$, in such a manner as to further reduce the error.  As shown in $(9)$, in this second extrapolation, $\Omega_{10,2}$ is simply a combination of previous extrapolation factors such that:

\begin{align*}
\Omega_{10,2} = \Omega_{10,5}\Omega_{5,2} = \left(\dfrac{10}{5}\right)^2\left(\dfrac{5}{2}\right)^2=\left(\dfrac{10}{2}\right)^2.
\end{align*}

\noindent Therefore, our next extrapolation will be $A_{10,2}$ as computed below:

\begin{align*}
A_{10,2} = \dfrac{\left(\dfrac{10}{2}\right)^2 12500833.\overline{3}-12520833.\overline{3}}{\left(\dfrac{10}{2}\right)^2-1} = 12500000.
\end{align*}

\noindent This estimate is exact because we have reduced the error to $O(n^2h^8)$.  This error term contains the difference of the seventh derivative of the function evaluated at $a$ and the seventh derivative evaluated at $b$.  However, that derivative is just a constant, allowing $f^{(7)}(a) = f^{(7)}(b)$, which then allows this error term to negate itself.\\
\\
\\
Example 2:\\

\noindent This example was chosen because it is continuously differentiable, and thus, using the alternative integration method presented in this paper, we will always have some non-zero error term. It also allows for four separate composites to be computed, and we can more readily observe the role of $h$, because in this case, $h \neq 1$. \\

\noindent Consider $f(x) = sin(x)$.  Then $\int_\pi^{2\pi} f(x)\ dx = -2$.  Assume, however, that we do not know the function, but are instead given 13 equally spaced points from $x = \pi$ to $2\pi$ and their corresponding values for $f(x)$, rounded to 10 decimal points.

\begin{table}[H]
\centering
\scalebox{1}{
\begin{tabular}{l l l}
$x$ & & $f(x)$ \\ [0.5ex]
\hline
$\pi$ & & 0\\
$13\pi/12$ & & -0.2588190451 \\
$7\pi/6$ & & -0.5 \\
$15\pi/12$ & & -0.7071067812 \\
$4\pi/3$ & & -0.8660254038 \\
$17\pi/12$ & & -0.9659258263 \\
$3\pi/2$ & & 1 \\
$19\pi/12$ & & -0.9659258263 \\
$5\pi/3$ & & -0.8660254038 \\
$21\pi/12$ & & -0.7071067812 \\
$11\pi/6$ & & -0.5 \\
$23\pi/12$ & & -0.2588190451 \\
$2\pi$ & & 0 \\ [1ex]
\hline
\end{tabular}}
\label{table:2}
\end{table}

\noindent Again, we calculate $A_{12}$ by first calculating $U_\infty$ and $O_\infty$.  However, note that $n = 12$, and therefore, $h = (b-a)/n = \pi/12$.  Also, note that $H = hf(a) = 0$.  Then, as stated previously, we calculate $U_1$ as the underestimate, which is $h$ multiplied by the sum of the first 12 values of $f(x)$, which gives us $U_1 = -1.9885637766$.  Likewise, the overestimate becomes $O_1 = -1.9885637766$.  Then we have

\begin{align}
U_{\infty}  = -1.9885637766 + \dfrac{-1.9885637766}{11} - \dfrac{12}{11}\cdot 0 = -2.1693423017, \\
O_{\infty} = -1.9885637766 - \dfrac{-1.9885637766}{13} + \dfrac{12}{13}\cdot 0 = -1.8355973322.
\end{align}

\noindent Averaging $(16)$ and $(17)$ then yields the initial approximation, $A_{12} = -2.0024698170$.\\

\noindent Now we can construct the various available composites.  We have four options when $n = 12$, a composite constructed using two, three, four, or six parts of six, four, three, or two subintervals each, respectively.  Denoting each composite with the number of subintervals used in each part during its construction gives us $A_2$, $A_3$, $A_4$, and $A_6$.  We would expect $A_2$ to be the most accurate but most difficult to compute, as it is a composite of many parts, and $A_6$ to be the least accurate but easiest to compute.\\

\noindent Following the same procedure as in example 1, $A_6$ is calculated as follows:\\

\noindent The first underestimate is found using values f($\pi$) through f($17\pi/12$), and our first overestimate is found using values f($13\pi/12$) through f($3\pi/2$).  Then our second underestimate is found using values f($3\pi/2$) through f($23\pi/12$), and our second overestimate is found by using values f($19\pi/12$) through f($2\pi$).  Then we have

\begin{align}
U_{\infty} & = -0.8633821944 + \dfrac{-0.8633821944}{5} - \dfrac{6}{5}\cdot 0 = -1.0360586333, \\
O_{\infty} & = -1.1251815822 - \dfrac{-1.1251815822}{7} + \dfrac{6}{7}\cdot 0 = -0.9644413562.
\end{align}

\noindent Averaging $(18)$ and $(19)$ gives us $-1.0002499947$, which is our approximation of the area of the first half of the entire integral.  Note that, as in the previous example, we must now calculate a new $H$, as $f(a)$ in this second half differs from that in the first half.  Then we have $H_2 = hf(3\pi/2) = -0.2617993878$, and

\begin{align}
U_{\infty} & = -1.1251815822 + \dfrac{-1.1251815822}{5} - \dfrac{6}{5}\cdot -0.2617993878 = -1.0360586333,\\
O_{\infty} & = -0.8633821944 - \dfrac{-0.8633821944}{7} + \dfrac{6}{7}\cdot -0.2617993878 = -0.9644413562.
\end{align}

\noindent Averaging $(20)$ and $(21)$ gives us $-1.0002499947$, which is our approximation of the area of the second half of the entire integral.  Adding this to the first half approximation gives us $A_6 = -2.0004999894$.\\

\noindent Now that we have our initial approximation and a composite approximation, we can find the appropriate extrapolation factor, $\Omega_{12,6}$.  Again, as this is our first extrapolation, we know that the numerator of $\Omega_{12,6}$ is the total number of subintervals in the original approximation and that the denominator is the number of subintervals used in each iteration of the composite.  Therefore, in this case, $\Omega_{12,6} = (n/m)^2=(12/6)^2$.\\

\noindent We calculate our new approximation, $A_{12,6}$, as follows:

\begin{align*}
A_{12,6} = \dfrac{\left(\dfrac{12}{6}\right)^2A_6 - A_{12}}{\left(\dfrac{12}{6}\right)^2-1} =\dfrac{\left(\dfrac{12}{6}\right)^2(-2.0004999894) - (-2.0024698170)}{\left(\dfrac{12}{6}\right)^2-1} = -1.9998433802.
\end{align*}

\noindent Using the same procedure to calculate the other available composites gives us $A_2 = -2.0000526243$, $A_3 = -2.0001193864$, and $A_4 = -2.0002147374.$  Now that we have our composite estimates, all that is required is to follow the previously described method to calculate our $\Omega$ and our extrapolations, and we can construct Table 2 similar to how Table 1 was constructed in the procedure section.

\begin{table}[H]
\caption{}
\centering
\scalebox{.8}{
\begin{tabular}{l r r r r r}
Error & $O(n^2h^4)$ & $O(n^2h^6)$ & $O(n^2h^8)$ & $O(n^2h^{10})$ & $O(n^2h^{12})$ \\[0.5ex]
\hline
$A_{12}$ & -2.0024698170 & -1.9998433802 & -2.0000010844 & -1.9999999828 & -2.0000000005 \\

$A_6$ & -2.0004999894 & -1.9999967037 & -2.0000000517 & -1.9999999985 & \\

$A_2$ & -2.0000526243 & -1.9999992147 & -2.0000000221 & \\

$A_3$ & -2.0001193864 & -1.9999967923 & \\

$A_4$ & -2.0002147374 & \\ [1ex]
\hline
\end{tabular}}
\label{table:nonlin2}
\end{table}

\noindent Note that in both examples, $f(a)$ was not the lowest value for $f(x)$ over the entire integral.  The terms ``underestimate'' and ``overestimate'' are, in some cases, misnomers.  If the function is strictly decreasing over the integral, nothing has to be changed in the calculation of $A_n$ or $A_m$.  This would result in the underestimate being larger than the overestimate, but this should not be alarming because the results of this procedure are not affected.

\section{Comparison with Romberg Integration}
\noindent Let us consider the results of Example 1.  If we had decided instead to use the Romberg method of integration and were unable to increase the number of points, we would have had an initial estimate of 50,000,000 using only the endpoints and the trapezoidal rule.  Our second estimate would have used the center point, essentially constructing a composite trapezoidal rule of two parts, but would not have fared much better at 25,390,625.  If we had then used the extrapolation method inherent in Romberg integration, our estimate would have been reduced to 17,187,500, but we would be unable to continue further at this point, with our best estimate having an error of almost 4.7 million. \\

\noindent This is perhaps not a fair assessment of the Romberg method, because the polynomial used is of high degree, which can result in large values for the $O(h^4)$ error term.  Furthermore, the Romberg method is crippled by being unable to use the optimal number of points.  Then let's allow the number of points to increase to 17 such that $h=10/16$.  Counting the original trapezoidal estimation, Romberg integration will be able to use five separate approximations.  In this case, Romberg integration will be able to easily reach $O(h^8)$ error to arrive at an exact answer.  \\

\noindent In fact, the Romberg method would also be able to arrive at this accuracy with only 9 points, or 8 subintervals.  However, with those same 8 subintervals, the alternative integration method presented in this paper would also arrive at the exact value.  Unfortunately though, Romberg integration cannot reach this same accuracy with any number of subintervals less than 8, or between 8 and 16.  It is only able to effectively use composites constructed of parts where the number of subintervals used in each part are powers of 2.  Considering only the cases where we are given 16 or fewer subintervals, our alternative method would be able to reach an accuracy of $O(n^2h^8)$.  In this example, we would arrive at the exact answer, given 6, 8, 10, 12, 14, 15 or 16 subintervals, with 12 and 16 subintervals providing even higher degrees of accuracy.\\

\noindent Let us then consider Example 2.  With the given 13 points, or 12 subintervals, Romberg integration will only be able to arrive at an estimate of -1.998570731824, which is a result of $O(h^6)$ error.  It would make use of the approximations from the normal trapezoidal rule, the composite trapezoidal rule composed of two parts, and the composite trapezoidal rule composed of four parts.  Any further composites that could be used are not available.\\

\noindent Suppose, again, that this is not a fair assessment.  Let's provide the Romberg method with an ideal number of subintervals, 32.  In this case, the Romberg method will be able to reach the error term $O(h^{12})$.  As is commonly done when implementing Romberg integration, we could build the following table of results.

\begin{table}[H]
\caption{}
\centering
\scalebox{.7}{
\begin{tabular}{l r r r r r r}
$h$ & $O(h^2)$ & $O(h^4)$ & $O(h^6)$ & $O(h^8)$ & $O(h^{10})$ & $O(h^{12})$ \\[0.5ex]
\hline
$\pi$ & 0 & -2.09439510239320 & -1.99857073182384 & -2.00000554997968 & -1.99999999458730 & -2.00000000000133 \\
$\pi/2$ & -1.57079632679490 & -2.00455975498443 & -1.99998313094599 & -2.00000001628805 & -1.99999999999604 \\
$\pi/4$ & -1.89611889793705 & -2.00026916994839 & -1.99999975245458 & -2.00000000005968 \\
$\pi/8$ & -1.97423160194556 & -2.00001659104794 & -1.99999999619085 \\
$\pi/16$ & -1.99357034377234 & -2.00000103336942 \\
$\pi/32$ & -1.99839336097015 \\ [1ex]
\hline
\end{tabular}}
\label{table:nonlin2}
\end{table}

\noindent The alternative method introduced in this paper can use the same composites, except for the one using 32 subintervals, as that composite is incorporated in the original approximation, $A_{32}$.  However, because the Romberg method's initial estimations have error of $O(h^2)$, and the method introduced in this paper has initial estimations with an error of $O(n^2h^4)$, these two methods are easily comparable.  Similarly as to how Table 3 was created, using the integration method presented in this paper, we can create Table 4, which, in practice, may be useful in keeping track of the various values of our extrapolation factors and approximations.

\begin{table}[H]
\caption{}
\centering
\scalebox{.8}{
\begin{tabular}{l r r r r r}
Error & $O(n^2h^4)$ & $O(n^2h^6)$ & $O(n^2h^8)$ & $O(n^2h^{10})$ & $O(n^2h^{12})$ \\[0.5ex]
\hline
$A_{32}$ & -2.00034682466611 & -1.99999967732512 & -2.00000000125221 & -1.99999999998018 & -2.00000000000133 \\

$A_2$ & -2.00000103336942 & -1.99999999619085 & -2.00000000005968 & -1.99999999999604 & \\

$A_4$ & -2.00000414490512 & -1.99999993815840 & -2.00000000406904 & \\

$A_8$ & -2.00001676514528 & -1.99999894949885 & \\

$A_{16}$ & -2.00007021208456 & \\ [1ex]
\hline
\end{tabular}}
\label{table:nonlin2}
\end{table}

\noindent Though the final approximations reached by both methods is exactly the same, this does not appear to always be the case in every similar situation.  However, when the estimation methods provide different results using the same composites, it is currently unclear as to what degree the differences, often at the 15th or further digit, could be due to round-off error.  Additionally, we were able to reach this final approximation as the 15th approximation, where the Romberg method uses 21 approximations to arrive at the same value. \\

\noindent Furthermore, note that if we had any number of subintervals less than 32, we would not be able to achieve this level of accuracy with Romberg integration.  In fact, we would not be able to achieve this level of accuracy again until we were given 64 subintervals.  However, our alternative method would be able to reach this level of accuracy given 12, 18, 20, 24, 28, or 30 subintervals, for those less than 32, and even an accuracy of $O(n^2h^{16})$ given either 24 or 30 subintervals.  A thorough discussion of the principles of Romberg integration can be found in [3].

\section{Conclusion}
\noindent The integration method introduced in this paper provides a manner by which to approximate the definite integral to a high degree of accuracy without it becoming necessary to include additional information.  The Romberg method of integration, though powerful when given a number of subintervals that is a power of two, is severely weakened otherwise.  Once the accuracy limit is reached, the standard approach in the Romberg method is to double the given number of points, thereby providing an additional composite and thus an additional level of accuracy.  \\

\noindent The Romberg method of using composites that increasingly double the number of subintervals used in their construction, for each consecutive extrapolation, suggests that the alternative method introduced in this paper may be at least as accurate as Romberg's method in situations optimal for the Romberg method, but many times more accurate otherwise, since our alternative method provides the added flexibility of using data from other available composites.  The Romberg method is particularly crippled when the number of subintervals available is odd.  However, our alternative method would be successful, as composites may still be constructed in these circumstances.  It appears that our alternative method provides a higher level of accuracy than the Romberg method when used on any number of subintervals that is not a power of two, and it may at least match the Romberg method's accuracy otherwise, and with fewer necessary approximations.\\

\noindent As mentioned previously, though only four composites were used in the most thorough example in this paper, there is no reason to believe that this method is restricted to an error no smaller than $O(n^2h^{12})$.  Providing that an additional unused composite can be computed, this error term should be able to be eliminated as well.  The accuracy of this method appears to be dependent on the number of composites available.

\section{Acknowledgements}
\noindent I would like to thank Professor Adam Glesser, PhD, in the Department of Mathematics at California State University, Fullerton, who shared his time and patience with me as he helped to confirm the validity of my method as introduced in this paper.  He is extremely dedicated, and his enthusiasm, support, and guidance were instrumental in providing the motivation necessary to complete this research.
 
\section{References}
\noindent$[1]$ Pedro Americo Almeida Magalhaes Junior and Cristina Almedia Magalhaes.  Higher-Order Netwon-Cotes Formulas.  Journal of Mathematics and Statistics, 6(2):193-204, 2010.\\
$[2]$  Konrad Knopp.  Theory and Application of Infinite Series.  Blackie and Son Limited, London, 1951. \\ 
$[3]$ Ole Østerby.  Romberg Integration:  Extrapolation and Error Estimation.  Aarhus University, 2005.\\
$[4]$ Avram Sidi.  Practical Extrapolation Methods.  Cambridge University Press, UK, 2003.\\
$[5]$  Charles A. Thompson. A Study of Numerical Integration Techniques for use in the companion circuit method of transient circuit analysis.  Purdue University ECE Technical Reports, Paper 297, 1997.\\
$[6]$  Thomas W. Tucker.  Rethinking Rigor in Calculus: The Role of the Mean Value Theorem.  Amer. Math. Monthly 104 (1997), 231-240. \\

\end{document}